\documentclass{artjlt}

\usepackage{amsmath,amsfonts,amssymb}
\usepackage{latexsym, graphicx, enumerate, amsthm}
\usepackage[all]{xy}

\newcommand{\Fn}{\mathfrak{n}}

\newcommand{\Fv}{\mathfrak{v}}
\newcommand{\Fz}{\mathfrak{z}}
\newcommand{\Ft}{\mathfrak{t}}

\newcommand{\RR}{\mathbb{R}}
\newcommand{\ZZ}{\mathbb{Z}}

\newcommand{\GG}{\mathcal{G}}

\theoremstyle{definition}
\newtheorem{constr}[Theorem]{Construction}

\title{Two-Step and Three-Step Nilpotent Lie Algebras Constructed from Schreier Graphs}                                     
\author{Allie Ray}                 
\lastname{Ray}  

\msc{05C99, 17B30, 22E25}    

\keywords{Metric Nilpotent Lie Algebras, Schreier Graphs, Gassmann-Sunada Triples}         

\address{Allie Ray\\University of Texas at Arlington\\Math Dept Box 19408\\Arlington, TX 76019-0408\\            
allieray@uta.edu              
}

\begin{document}


\maketitle

\begin{abstract}
We  associate a two-step nilpotent Lie algebra to an arbitrary Schreier graph. We then use properties of the Schreier graph to determine necessary and sufficient conditions for this Lie algebra to extend to a three-step nilpotent Lie algebra. As an application, if we start with pairs of non-isomorphic Schreier graphs coming from Gassmann-Sunada triples, we prove that the pair of associated two-step nilpotent Lie algebras are always isometric. In contrast, we use a well-known pair of Schreier graphs to show that the associated three-step nilpotent extensions need not be isometric.   
\end{abstract}

\section{Introduction} \label{intro}
The purpose of this paper is to introduce a new method for associating two-step nilpotent Lie algebras with Schreier graphs and then extend these Lie algebras (in certain cases) to three-step nilpotent Lie algebras. Also, we compare the pairs of resultant Lie algebras when associated with Gassmann-Sunada triples.\\

In 2004, S.G. Dani and M.G. Mainkar first presented a method for constructing two-step nilpotent Lie algebras from simple graphs \cite{DM}. They used the two-step nilpotent  construction to find properties of a graph that would result in the constructed manifold admitting Anosov automorphisms. J. Lauret and C. Will used this construction to find examples of nonisometric Einstein solvmanifolds \cite{LW}, and H. Pouseele and P. Tirao used the construction to consider symplectic nilmanifolds \cite{PT}. Mainkar also proved that for simple graphs, the resulting Lie algebras are isomorphic if and only if the graphs are isomorphic \cite{MM}, and in \cite{kstep}, she extended this construction to $k$-step nilpotent Lie algebras. Also, V. Grantcharov is currently working on extending the Dani-Mainkar construction on simple graphs to three-step solvable Lie algebras \cite{grant}. In the Dani-Mainkar construction, each vertex and each edge of the graph correspond to distinct elements in the Lie algebra; therefore for large graphs, the corresponding dimension of the Lie algebra is also large. For the higher-step construction, the dimension of the constructed Lie algebras grows more rapidly.  \\

The focus of this paper is on Schreier graphs because of their inherent group structure. J.L. Gross proved that every connected regular graph of even degree is a Schreier graph, \cite{gross}. Schreier graphs; however, are often non-simple directed graphs, in which case the Dani-Mainkar construction is not defined. We therefore introduce a new method for associating Lie algebras with Schreier graphs, as suggested by C.S. Gordon. In \S\ref{defns}, we discuss the definitions and notations that will be used in this paper. In \S\ref{23step}, we detail this new construction of a two-step nilpotent Lie algebra associated with an arbitrary Schreier graph. We also provide necessary and sufficient conditions on the graph for this construction to extend to a three-step nilpotent Lie algebra. As an application, in \S\ref{lags} we prove that for any pair of Schreier graphs associated to a Gassmann-Sunada triple, the resultant two-step nilpotent Lie algebras are isometric, but we give an example where the pair of three-step nilpotent Lie algebras are non-isometric.\\


\section{Background Info and Notation} \label{defns}
\begin{Definition}\label{schr} Let $G$ be a finite group and $H$ a subgroup of $G$. Let $C:=\{z_1,\dots, z_c,z_1^{-1},\dots,z_c^{-1}\}$ be a generating set of $G$ not containing the identity that is closed under inverses, and let $C_{pos}:=\{z_1,\dots,z_c\}$. The \emph{Schreier graph} of $G$ relative to $H$ and $C$, written $\GG(G,H,C)$ or simply $\GG$ if understood in context, is a directed edge-labeled graph defined by the following. The vertices of $\GG$ consist of the set of right cosets, $V(\GG)=\{Hg:\; g \in G\}.$ The edges consist of the set of ordered pairs $E(\GG)=\{(H g,H g z_i^{-1}): z_i \in C_{pos}\}$, and each edge $(Hg,Hg z_i^{-1})$ is given the label $z_i$. \end{Definition} 

\begin{Example} Let $G=S_4, \, H=S_3, \mbox{ and } C_{pos}=\{(123),\,(1234)\}$. Then the Schreier graph $\GG$ with respect to $H$ and $C$ is
$$\UseTips 
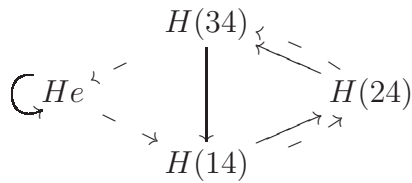
\xymatrix@R=.3cm{
       & H(34) \ar@//@{-->}[dl] \ar[dd] & \\
			He \ar@(ul,dl) \ar@//@{-->}[dr] & & H(24) \ar[ul] \ar@/_/@{-->}[ul]  \\
			&  H(14) \ar[ur] \ar@/_/@{-->}[ur] & }$$
where the solid lines correspond to edges formed by the first generator, $(123)$, in $C_{pos}$ and dotted lines to the second generator, $(1234)$.
\end{Example}

\begin{Remark}\label{zz-1} Note that while a Schreier graph is defined for an element of $C_{pos}$ of order 2, the edges associated to those elements will become trivial elements in the Lie algebras we construct in \S\ref{23step}. Hence, in what follows, we assume our generating set $C$ does not contain order 2 elements, i.e. $z \neq z^{-1} \mbox{ for all } z \in C$. \end{Remark}
\begin{Remark}\label{action} The structure of a Schreier graph implies that the group $G$ acts on $V(\GG)$ by right inverse multiplication. To see this, we define $\alpha(z_i): V(\GG) \mapsto V(\GG)$ for $z_i \in C$ by $$\alpha(z_i)(Hg)=Hgz_i^{-1} \mbox{ for all } z_i \in C.$$ Then $\alpha$ extends from $C$ to $G$ because $C$ generates $G$. \end{Remark}

\begin{Definition} For a directed graph, a \emph{walk} of length $q$ from vertex $v$ to $w$ is a sequence of $q+1$ vertices (and therefore $q$ edges) where successive vertices in the sequence are connected by a directed edge in $E(\GG)$. If these vertices are all distinct, except possibly $v$ and $w$, then this is called a \emph{path} of length $q$, or a \emph{q-path}. If $v=w$, then this path is called \emph{closed}. We will denote $q$-paths by a $(q+1)$-tuple of vertices, $(v_1, v_2, \dots, v_{q+1})$ where $(v_i,v_{i+1})\in E(\GG)$ for $i=1,\dots,q$. \end{Definition}

The following two facts about properties of a Schreier graph will be important in the proofs of the main theorems in \S\ref{23step} and \S\ref{lags}.
\begin{Remark}\label{closedpaths} Because the edges of a Schreier graph are associated with generators of a finite group, Schreier graphs will be the union of closed paths of a single label, where the length of a path with every edge labeled $z_i \in C_{pos}$ is less than or equal to the order of the generating element $z_i$. \end{Remark}
\begin{Remark}\label{possibilities} If $|C_{pos}|=c$, then the Schreier graph $\GG$ will be $2c$-regular where each vertex has a directed edge labeled $z_i$ going out of the vertex and one going into the vertex, $i=1,\dots,c$. This gives three different possibilities for each vertex $v$ and each generator (and hence each label) $z$:\\ 
%
\begin{tabular}{ccc} 1. $\alpha(z)(v) \neq \alpha(z^{-1})(v)$ & 2. $\alpha(z)(v) = \alpha(z^{-1})(v) \neq v$, & and 3. $\alpha(z)(v) = \alpha(z^{-1})(v)= v$\\
$\UseTips 
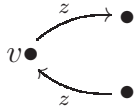
\xymatrix@R=.1cm{ & \bullet\\
       v\bullet \ar@/^/[ur]^z & \\
        & \bullet  \ar@/^/[ul]^z}$
& $\UseTips 
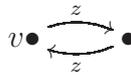
\xymatrix@R=.1cm{ & \\
       v\bullet \ar@/^/[r]^z & \bullet \ar@/^/[l]^z
        }$	
& $\UseTips 
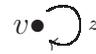
\xymatrix@R=.1cm{ & \\
       v \bullet \ar@(ur,dr)[]^z}$	
\end{tabular}
\end{Remark}

The following two definitions for Lie algebras are standard, see e.g. P. Eberlein \cite{eber}.
\begin{Definition} A Lie algebra $\Fn$  is $k$-step \emph{nilpotent} if $\Fn^{(k)}=0$ for some $k \in \ZZ^+$ but $\Fn^{(k-1)}\neq0$, where $\Fn^{(0)}=\Fn$ and $\Fn^{(k)}=[\Fn^{(k-1)},\Fn],\,k \in \ZZ^+$. \end{Definition}

\begin{Definition} Given two Lie algebras, $(\Fn_1, \, [\, , \,]_1)$ and $(\Fn_2, \, [\, , \,]_2)$, a map $\phi:\Fn_1 \to \Fn_2$ is an \emph{isomorphism} if $\phi$ is a linear bijection and $\phi([x,y]_1)=[\phi(x),\phi(y)]_2$ for all $x,y \in \Fn_1$. If an inner product is specified on these Lie algebras we say that $\phi$ is an \emph{isometry} if it is an isomorphism and in addition $<x,y>_1=<\phi(x),\phi(y)>_2$ for all $x,y \in \Fn_1$. \end{Definition}

In \S\ref{23step}, we develop a construction of two-step and three-step nilpotent Lie algebras associated with Schreier graphs, and in \S\ref{lags}, we compare properties of pairs of Lie algebras arising from Gassmann-Sunada triples. 

\begin{Definition} Let $G$ be a finite group, with $H_1$ and $H_2$ subgroups of $G$ such that for every $g \in G$, \begin{eqnarray} |[g] \cap H_1| = |[g] \cap H_2|, \end{eqnarray} where $[g]$ denotes the conjugacy class of $g$ in $G$. In this case, $H_1$ and $H_2$ are called \emph{almost conjugate} subgroups of $G$, and $(G,H_1,H_2)$ is called a \emph{Gassmann-Sunada triple}, \cite{Bus, Sun}. \end{Definition}

\begin{Example}\label{brookspic} Let $G=SL(3,2),\; H_1=\left\{ \scriptsize{ \left( \begin{array}{ccc}  1 & * & * \\ 0 & * & * \\ 0 & * & *  \end{array} \right)} \right\},\; \mbox{ and } H_2= \left\{ \scriptsize{ \left( \begin{array}{ccc} 1&0&0\\ *&*&*\\ *&*&* \end{array} \right)}   \right\}$.\\
Let $C_{pos}=\left\{ \scriptsize{\left( \begin{array}{ccc} 0&1&1\\0&1&0\\1&0&0 \end{array} \right),\; \left( \begin{array}{ccc} 1&0&0\\0&0&1\\0&1&1 \end{array} \right)} \right\}$. In \cite{Brooks, Bus}, this is shown to be a Gassmann-Sunada triple. The Schreier graphs, $\GG_1=\GG(G,H_1,C)$ and $\GG_2=\GG(G,H_2,C)$, are\\ 
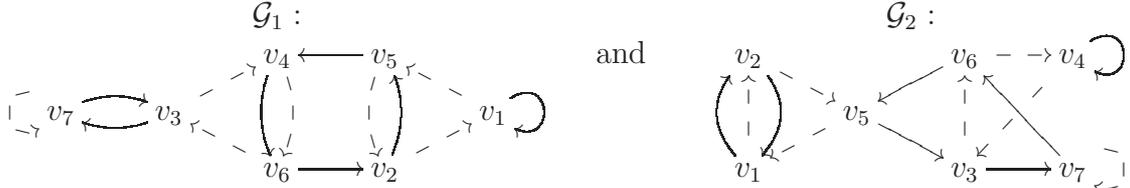
\begin{tabular}{ccc} $\GG_1:$&  &$\GG_2:$ \\
$\UseTips \xymatrix@R=.3cm{
       & & v_4 \ar@/_/[dd] \ar@/^/@{-->}[dd]   &v_5 \ar[l] \ar@/_/@{-->}[dd]  & \\
			v_7 \ar@//@{-->}@(ul,dl) \ar@/^/[r] &v_3 \ar@//@{-->}[ur] \ar@/^/[l]  & & &v_1\ar@//@{-->}[ul]  \ar@(ur,dr) \\
			& & v_6 \ar@//@{-->}[ul] \ar[r] &v_2  \ar@/_/[uu] \ar@//@{-->}[ur] &}$  
& \hspace{.5cm} and \hspace{.5cm} &
$\UseTips \xymatrix@R=.3cm{
v_2 \ar@//@{-->}[dr] \ar@/^1pc/[dd] & & v_6 \ar@//@{-->}[r] \ar[dl] & v_4 \ar@//@{-->}[ddl] \ar@(ur,dr)\\
& v_5 \ar@//@{-->}[dl] \ar[dr] & & \\
v_1 \ar@//@{-->}[uu] \ar@/^1pc/[uu]  & & v_3 \ar@//@{-->}[uu] \ar[r]  & v_7 \ar[uul]  \ar@//@{-->}@(ur,dr)}$
\end{tabular}\\
where the solid line corresponds to the first generator given in $C_{pos}$ and the dotted line to the second generator.
\end{Example}


\section{Two-Step and Three-Step Nilpotent Lie Algebra Constructions} \label{23step}

In this section, we construct a two-step nilpotent Lie algebra from an arbitrary Schreier graph and then give necessary and sufficient conditions under which this Lie algebra can extend to a three-step nilpotent Lie algebra.\\

\begin{constr}[Two-Step Nilpotent Construction]\label{constr2}
From a Schreier graph $\GG=\GG(G,H,C)$ given by Definition  \ref{schr}, we let $\Fv$ be the space of formal linear combinations over $\RR$ of elements in $V(\GG)$ and $\Fz$ be the space of formal linear combinations over $\RR$ of elements in $C_{pos}$, where $|C_{pos}|=c$. We then define the Lie algebra $\Fn:= \Fv \oplus \Fz$ as the direct sum of vector spaces; we then require $\Fz$ to be contained in the center of $\Fn$ and define the Lie bracket by the following: $\forall v_i,v_j \in V(\GG) \subseteq \Fv$,
\begin{align} [v_i,v_j] &= \sum_{p=1}^{|C_{pos}|} (\epsilon_p - \epsilon'_p) z_p, \label{2bracket}\\
\mbox{where } \epsilon_p &= \left\{ \begin{array}{rcl} 1 &,& \mbox{ if } v_j=\alpha(z_p)(v_i)\\ 0 &,& \mbox{ otherwise,} \end{array} \right. \notag\\
 \mbox{and } \epsilon'_p &= \left\{ \begin{array}{rcl} 1 &,& \mbox{ if } v_j=\alpha(z_p^{-1})(v_i)\\ 0 &,& \mbox{ otherwise.} \end{array} \right. \notag \end{align}
All other brackets not defined by linearity or skew-symmetry are set equal to zero. \end{constr}

To see that this does define a Lie algebra, consider the following. First, if $z=z^{-1}\in C_{pos}$, then $[v_i,v_j]=0 \; \forall v_i,v_j \in \Fv$, which is why we exclude such elements from $C_{pos}$ as we mentioned in Remark~\ref{zz-1}. Also note that $[v_i,v_i]=0$ because for a fixed label $z_p$, either $v_i$ has a loop with label $z_p$ in which case $\epsilon_p=\epsilon'_p=1$, or $v_i$ does not have a loop with label $z_p$ in which case $\epsilon_p=\epsilon'_p=0$. In either case, $\epsilon_p-\epsilon'_p=0$ for all $p$.  Furthermore, this bracket will be skew-symmetric because $v_j=\alpha(z_p)(v_i)$ implies $v_i=\alpha(z_p^{-1})(v_j)$. Finally, note that because $\Fz$ is contained in the center of $\Fn$, the Jacobi identity on the bracket given above is trivial, which makes $(\Fn,[\, , \,])$ as defined above a two-step nilpotent Lie algebra. 

\begin{Remark}\label{onb} In this paper when needed, we specify an inner product on $\Fn=\Fv \oplus \Fz$ by requiring $\{V(\GG),C_{pos}\}$ to be an orthonormal basis. \end{Remark}

\begin{Remark} The two-step nilpotent Lie algebra defined in Construction \ref{constr2} does not rely on the fact that the graph was a Schreier graph.  A two-step nilpotent Lie algebra can be constructed similarly from any directed, labeled (colored) graph by having a set of graph labels (colors), $C_{pos}=\{z_1,\dots,z_c\}$, instead of having a set of generators of a group acting on the graph. The Lie bracket on $\Fn:= \Fv \oplus \Fz$ is then defined as in Construction \ref{constr2}, except now
\begin{align} \epsilon_p &= \left\{ \begin{array}{rcl} 1 &,& \mbox{ if } (v_i,v_j) \mbox{ is an edge labeled } z_p\\ 0 &,& \mbox{ otherwise,} \end{array} \right. \notag\\
 \mbox{and } \epsilon'_p &= \left\{ \begin{array}{rcl} 1 &,& \mbox{ if } (v_j,v_i) \mbox{ is an edge labeled } z_p\\ 0 &,& \mbox{ otherwise.} \end{array} \right. \notag \end{align} \end{Remark}

In order to find necessary and sufficient conditions on a Schreier graph for the two-step nilpotent Lie algebra construction to extend to a three-step nilpotent Lie algebra, we must introduce the following definition.

\begin{Definition}\label{admissible} For a Schreier graph $\GG=\GG(G,H,C)$, a label $z \in C_{pos}$ is called \emph{admissible} if there exists a single closed path of length 3 or 4 with each edge labeled $z$, and all other closed paths with edges labeled $z$ are of length 1 or 2. Otherwise, $z$ is called \emph{inadmissible}. We will denote the set of admissible labels by $\{z_{r_1},\dots,z_{r_m}\}$ and the set of inadmissible labels by $\{z_{b_1},\dots,z_{b_n}\}$. A path is called \emph{admissible} if it is the single closed path of length 3 or 4 for an admissible label $z_r$.\end{Definition}

\begin{Theorem}\label{3stepthm} Let $G$ be a finite group, $H$ a subgroup of $G$, $C$ a generating set of $G$, and $\GG$ the Schreier graph of $G$ with respect to $H$ and $C$ as in Definition~\ref{schr}. Let $\Fn$ be the two-step nilpotent Lie algebra associated with $\GG$ by Construction \ref{constr2}. Then $\Fn$ extends to a three-step nilpotent Lie algebra $\widehat{\Fn}$ if and only if there exists at least one admissible label in $C_{pos}$. Moreover, up to the variations allowed in Construction \ref{constr3} below, this is the only 3-step nilpotent extension of $\Fn$. \end{Theorem}

\begin{constr}[Three-Step Nilpotent Construction]\label{constr3}
For each admissible label $z_{r_k}$, we define new elements $t_{r_{k,1}}$ and $t_{r_{k,2}}$ (at least one $t_{r_{k,\ell}}\neq 0$) such that the 3-step nilpotent extension of $\Fn$ is $\widehat{\Fn}=\Fv\oplus\Fz\oplus\Ft,$ where $\Fv$ and $\Fz$ are defined as before and $\Ft=\mbox{span}_{\RR} \{t_{r_{k,1}},\,t_{r_{k,2}}: \; z_{r_k} \mbox{ is admissible} \}$. The Lie bracket is then defined as in Construction \ref{constr2} with the following additional nonzero brackets, and then extend by linearity and skew-symmetry:\\
If the admissible path with label $z_{r_k}$ is of length 4 and has successive vertices $(v_1,v_2,v_3,v_4,v_1)$, we set 
\begin{align} \begin{array}{rcl} [v_1,z_{r_k}] = -[v_3,z_{r_k}] &=& t_{r_{k,1}}, \mbox{ and}\\
\left[v_2,z_{r_k}\right] = -[v_4,z_{r_k}] &=& t_{r_{k,2}} \end{array} \label{4loop} \end{align}
If the admissible path with label $z_{r_k}$ is of length 3 and has successive vertices $(v_1,v_2,v_3,v_1)$, we set
\begin{align} \begin{array}{rcl} [v_1,z_{r_k}] &=& t_{r_{k,1}},\\
\left[v_2,z_{r_k}\right] &=& t_{r_{k,2}}, \mbox{ and}\\
\left[v_3,z_{r_k}\right] &=& -(t_{r_{k,1}}+t_{r_{k,2}}) \end{array} \label{3loop} \end{align}
For any other vertex $v_i$ not in the admissible 3- or 4-path, we set
\begin{align} [v_i,z_{r_k}]=0 \label{vlzr} \end{align}
For any edge with inadmissible label $z_b$, we set
\begin{align} [v_j,z_b] &= 0 \; \forall v_j \in \Fv.\label{black} \end{align} \end{constr}

\begin{Remark} In order for $\widehat{\Fn}$ to be 3-step nilpotent, we must set at least one $t_{r_{k,\ell}} \neq 0$. The 3-step nilpotent extension of $\Fn$ is not unique.  Distinct Lie algebra extensions can be obtained by defining relations between the various elements $t_{r_{k,\ell}}$, namely these elements may be linearly dependent.  Because of these variations, we get $1 \leq \dim \Ft \leq 2m$, where $m$ is the number of admissible labels. \end{Remark}

\begin{Remark} This paper does not address extensions where $[v_i,v_j] \in \Ft$ since these do not seem to intuitively arise from graph properties, nor do they contribute to the extension being 3-step nilpotent. \end{Remark}

\begin{Example}\label{brooksex3} The following is a three-step nilpotent extension of the Lie algebras associated with the Schreier graphs in Example~\ref{brookspic}:\\
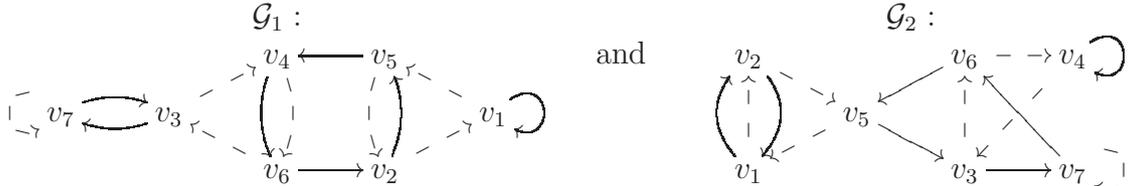
\begin{tabular}{ccc} $\GG_1:$&  &$\GG_2:$ \\
$\UseTips \xymatrix@R=.3cm{
       & & v_4 \ar@/_/[dd] \ar@/^/@{-->}[dd]   &v_5 \ar[l] \ar@/_/@{-->}[dd]  & \\
			v_7 \ar@//@{-->}@(ul,dl) \ar@/^/[r] &v_3 \ar@//@{-->}[ur] \ar@/^/[l]  & & &v_1\ar@//@{-->}[ul]  \ar@(ur,dr) \\
			& & v_6 \ar@//@{-->}[ul] \ar[r] &v_2  \ar@/_/[uu] \ar@//@{-->}[ur] &}$  
& \hspace{.5cm} and \hspace{.5cm} &
$\UseTips \xymatrix@R=.3cm{
v_2 \ar@//@{-->}[dr] \ar@/^1pc/[dd] & & v_6 \ar@//@{-->}[r] \ar[dl] & v_4 \ar@//@{-->}[ddl] \ar@(ur,dr)\\
& v_5 \ar@//@{-->}[dl] \ar[dr] & & \\
v_1 \ar@//@{-->}[uu] \ar@/^1pc/[uu]  & & v_3 \ar@//@{-->}[uu] \ar[r]  & v_7 \ar[uul]  \ar@//@{-->}@(ur,dr)}$
\end{tabular}\\
The solid lines correspond to the first generator in $C_{pos}$, denoted $z_r$ because it is admissible, and the dotted lines correspond to the second generator, denoted $z_b$ because it is inadmissible. If we delete the last column of bracket relations below, we have the two-step nilpotent Lie algebra as defined in Construction \ref{constr2}.\\
\begin{tabular}{r l l | l}
$\widehat{\Fn_1}: $ & $[v_1,v_2]=-z_b$ & $[v_3,v_4]=z_b$ & $[v_2,z_r]=t$\\
 & $[v_1,v_5]=z_b$ & $[v_3,v_6]=-z_b$ & $[v_4,z_r]=-t$\\
 & $[v_2,v_5]=z_r-z_b$ & $[v_4,v_5]=-z_r$ & $[v_5,z_r]=0$\\  
 & $[v_2,v_6]=-z_r$ & $[v_4,v_6]=z_r+z_b$ & $[v_6,z_r]=0$\\
 & & &\\
$\widehat{\Fn_2}: $ & $[v_1,v_2]=z_b$ & $[v_3,v_6]=z_b$ & $[v_3,z_r]=t$\\
 & $[v_1,v_5]=-z_b$ & $[v_3,v_7]=z_r$ & $[v_5,z_r]=0$\\
 & $[v_2,v_5]=z_b$ & $[v_4,v_6]=-z_b$ & $[v_6,z_r]=-t$\\  
 & $[v_3,v_4]=-z_b$ & $[v_5,v_6]=-z_r$ & $[v_7,z_r]=0$\\
 & $[v_3,v_5]=-z_r$ & $[v_6,v_7]=-z_r$ & \\ 
\end{tabular}\\
All other brackets not defined by skew-symmetry or linearity are equal to zero. \end{Example}

\begin{Proof}[Proof of Thm. \ref{3stepthm} (sufficient)]\\
Define $\epsilon_{i,j}^{r_k}=\left\{ \begin{array}{rcl} 1 &,& \mbox{ if there is a $z_{r_k}$-edge connecting $v_i$ to $v_j$}\\  -1 &,& \mbox{ if there is a $z_{r_k}$-edge connecting $v_j$ to $v_i$}\\ 0 &,& \mbox{ otherwise} \end{array} \right.$,\\ and similarly define  $\epsilon_{i,j}^{b_{\ell}}$. We proceed by induction on the number of admissible labels.
Assume that the Schreier graph has only one admissible label $z_r$, and the inadmissible labels, if any exist, are denoted $z_{b_{\ell}},\,\ell=1,\dots,n$. If we pick any three vertices from the graph, say $v_1,v_2,v_3$, then the following possibilities occur for the Jacobi identity on those three vertices:\\

Case 1: There are no edges labeled $z_r$ connecting $v_1,v_2$, or $v_3$, in which case the Jacobi identity will be satisfied because\\ 
$[v_1,[v_2,v_3]]+[v_2,[v_3,v_1]]+[v_3,[v_1,v_2]]$\\ 
$=[v_1,\epsilon_{2,3}^r z_r]+\sum_{\ell=1}^n [v_1,\epsilon_{2,3}^{b_{\ell}} z_{b_{\ell}}]+[v_2,\epsilon_{3,1}^r z_r]+\sum_{\ell=1}^n [v_2,\epsilon_{3,1}^{b_{\ell}} z_{b_{\ell}}]+[v_3,\epsilon_{1,2}^r z_r]+\sum_{\ell=1}^n [v_3,\epsilon_{1,2}^{b_{\ell}} z_{b_{\ell}}]$\\
\hspace*{\fill} by linearity of the bracket\\
 $=[v_1,0]+0+[v_2,0]+0+[v_3,0]+0$ \hfill by Equation~\ref{black} and definition of $\epsilon_{i,j}^r$\\
$=0$.\\
Note that by the linearity of the Lie bracket, we can always take the Jacobi identity and separate the brackets containing $z_{b_{\ell}}$ terms, which will equal zero by Equation~\ref{black}, so we only need to consider the Jacobi identity in relation to brackets containing $z_{r_k}$ terms.\\

Case 2: Without loss of generality, there is precisely one $z_r$-edge connecting $v_1$ to $v_2$, which implies that $v_3$ is not contained in the admissible path with label $z_r$.  In this case,\\
$[v_1,[v_2,v_3]]+[v_2,[v_3,v_1]]+[v_3,[v_1,v_2]]$\\
$=[v_1,0]+[v_2,0]+[v_3,z_r]$ \hfill by Equation \ref{black} and definition of $\epsilon_{i,j}^r$\\
$=0$ \hfill by Equation \ref{vlzr}.\\

Case 3: There are precisely two edges labeled $z_r$. Without loss of generality, one edge connects $v_1$ to $v_2$ and the other from $v_2$ to $v_3$. Since there is no $z_r$-edge connecting $v_3$ to $v_1$, this implies that the path with labels $z_r$ must be an admissible 4-path so $[v_1,z_r]=-[v_3,z_r]$ by Equation \ref{4loop}.  Again the Jacobi identity is satisfied because\\
$[v_1,[v_2,v_3]]+[v_2,[v_3,v_1]]+[v_3,[v_1,v_2]]$\\
$=[v_1,z_r]+[v_2,0]+[v_3,z_r]$ \hfill by Equation \ref{black} and definition of $\epsilon_{i,j}^r$\\
$=0$ \hfill by Equation \ref{4loop}.\\

Case 4: There are three edges labeled $z_r$ connecting $v_1$ to $v_2$ to $v_3$ back to $v_1$. So the path here is an admissible 3-path with label $z_r$.  The Jacobi equation becomes\\
$[v_1,[v_2,v_3]]+[v_2,[v_3,v_1]]+[v_3,[v_1,v_2]]$\\
$=[v_1,z_r]+[v_2,z_r]+[v_3,z_r]$ \hfill by Equation \ref{black} and definition of $\epsilon_{i,j}^r$\\
$=0$ \hfill by Equation~\ref{3loop}.\\
These four cases cover all possibilities because of the properties of a Schreier graph discussed in Remark~\ref{possibilities}. Therefore, no matter which three vertices we pick in the graph and by the linearity of the Lie bracket, the Jacobi identity is always satisfied, making $\widehat{\Fn}$ a Lie algebra.

Now using induction, assume that we have a Lie algebra associated with a graph with admissible labels, $z_{r_1},\dots,z_{r_m}$, and inadmissible labels, $z_{b_1},\dots,z_{b_n}$. If we add an additional admissible label $z_{r_{m+1}}$ in $C_{pos}$, then the Jacobi identity for any three vertices becomes\\
$[v_1,[v_2,v_3]]+[v_2,[v_3,v_1]]+[v_3,[v_1,v_2]]$\\
$=[v_1,\epsilon_{2,3}^{r_{m+1}} z_{r_{m+1}}]+\sum_{k=0}^m [v_1,\epsilon_{2,3}^{r_k} z_{r_k}]+[v_2,\epsilon_{3,1}^{r_{m+1}} z_{r_{m+1}}]+\sum_{k=0}^m [v_2,\epsilon_{3,1}^{r_k} z_{r_k}]+[v_3,\epsilon_{1,2}^{r_{m+1}} z_{r_{m+1}}]+\sum_{k=0}^m [v_3,\epsilon_{1,2}^{r_k} z_{r_k}]$ \hfill by Equation \ref{black} and linearity of the bracket\\
$=([v_1,\epsilon_{2,3}^{r_{m+1}} z_{r_{m+1}}]+[v_2,\epsilon_{3,1}^{r_{m+1}} z_{r_{m+1}}]+[v_3,\epsilon_{1,2}^{r_{m+1}} z_{r_{m+1}}])+\sum_{k=0}^m ([v_1,\epsilon_{2,3}^{r_k} z_{r_k}]+[v_2,\epsilon_{3,1}^{r_k} z_{r_k}]+[v_3,\epsilon_{1,2}^{r_k} z_{r_k}])$\\
$=[v_1,\epsilon_{2,3}^{r_{m+1}} z_{r_{m+1}}]+[v_2,\epsilon_{3,1}^{r_{m+1}} z_{r_{m+1}}]+[v_3,\epsilon_{1,2}^{r_{m+1}} z_{r_{m+1}}]+0$ \hfill by induction hypothesis.\\
$=0$ \hfill because the proof of the base case of the induction proof showed that the Jacobi identity is satisfied for any single admissible label. \end{Proof}

\begin{Proof}[Proof of Thm. \ref{3stepthm} (necessary)] Assume now that the Schreier graph $\GG$ has no admissible labels in $C_{pos}$. This means that for each label $z_{b_{\ell}}$, at least one of the following occur:\\
\begin{enumerate}
\item Each closed path with edges labeled $z_{b_{\ell}}$ is of length 1 or 2.
\item There are at least two closed paths of length 3 or 4, with edges labeled $z_{b_{\ell}}$.
\item There exists a closed path with label $z_{b_{\ell}}$ that is of length $q$, $q \geq 5$.
\end{enumerate}
We continue by induction on the number of inadmissible labels $z_{b_{\ell}}$ in the Schreier graph. Assume that $\GG$ only has one inadmissible label $z_b$.\\

Case 1: If each closed path in $\GG$ with label $z_b$ is of length 1 or 2, then $[v_i,v_j]=0$ for all $v_i,v_j \in \Fv$ by how the Lie bracket is defined in Equation \ref{2bracket}. Therefore, $\dim \Fz=0 \implies \dim \Ft=0$ so there does not exist a three-step nilpotent extension of $\Fn$.\\

Case 2: Assume that $\GG$ has at least two closed paths of length 3 or 4, with edges labeled $z_b$. Let $v_i$ be a vertex in one of these paths and $(v_j,v_k)$ be an edge in one of the other paths of length 3 or 4. Note that these two paths will not have any vertices in common by Remark~\ref{possibilities}. Because we are assuming that $\Fn$ is a Lie algebra and only considering when there is a 3-step nilpotent extension, we may assume that the Jacobi identity is satisfied for all $v \in \Fv$.  Therefore,
\begin{align} [v_i,[v_j,v_k]]+[v_j,[v_k,v_i]]+[v_k,[v_i,v_j]] &=0 \notag\\
\implies [v_i,z_b]+[v_j,0]+[v_k,0] &=0 \mbox{ (by Equation \ref{2bracket})} \notag\\
\implies [v_i,z_b] &= 0 \mbox{ for all $v_i$ in the closed path}. \label{case2a} \end{align}
Since this was for an arbitrary $v_i$ in a path of length 3 or 4, we can conclude that $[v_i,z_b]=0$ for all $v_i$ in any path of length 3 or 4. Now, let $v_i$ be another vertex on this graph not contained in a closed path of length 3 or 4, and again let $(v_j,v_k)$ be an edge in one of the closed paths of length 3 or 4. Then,
\begin{align} [v_i,[v_j,v_k]]+[v_j,[v_k,v_i]]+[v_k,[v_i,v_j]] &=0 \notag\\
\implies [v_i,z_b]+[v_j,0]+[v_k,0] &=0 \mbox{ (by Equation \ref{2bracket})} \notag\\
\implies [v_i,z_b] &= 0 \mbox{ for all $v_i$ not in the closed path of length 3 or 4}. \label{case2b} \end{align}
Therefore, $[v_i,z_b]=0$ for all $v_i\in \Fv$ (by Equations \ref{case2a} and \ref{case2b}), which implies that $\dim \Ft=0$ so a three-step extension of $\Fn$ of the type assumed does not exist.\\

Case 3: Assume that $\GG$ has a closed path of length $q, \; q \geq 5,$ with edges labeled $z_b$. Let the successive vertices of this closed path be $(v_0,v_1,\dots,v_{q-1},v_0)$. Because $\Fn$ is a Lie algebra, we assume that the Jacobi identity is satisfied for $v_i, v_{(i+2)\!\!\!\mod q}, \mbox{ and } v_{(i+3)\!\!\!\mod q}$. This implies that $\forall i=0,\dots,q-1$,\\
$[v_i,[v_{(i+2)\!\!\!\mod q},v_{(i+3) \!\!\!\mod q}]]+[v_{(i+2) \!\!\!\mod q},[v_{(i+3) \!\!\!\mod q},v_i]]+[v_{(i+3) \!\!\!\mod q},[v_i,v_{(i+2) \!\!\!\mod q}]] = 0$\\
$\implies [v_i,z_b]+[v_{(i+2) \!\!\!\mod q},0]+[v_{(i+3) \!\!\!\mod q},0] =0$ \hfill because two nonconsecutive points in a closed path with labels $z_b$ on a Schreier graph cannot have a $z_b$-edge connecting them.
\begin{align} \implies [v_i,z_b] =0 \; \forall i=0,\dots,q-1. \label{case3a} \end{align}
Now let $v_j$ be a vertex not in this closed path of length $q$. Then 
\begin{align} [v_j,[v_0,v_1]]+[v_0,[v_1,v_j]]+[v_1,[v_j,v_0]] &=0 \notag\\
\implies [v_j,z_b]+[v_0,0]+[v_1,0] &=0 \mbox{ (by Equation \ref{2bracket})} \notag\\
\implies [v_j,z_b] &= 0 \mbox{ for all $v_j$ not in the closed path of length $q$}. \label{case3b} \end{align}
Therefore, $[v_j,z_b]=0$ for all $v_j \in \Fv$ (by Equations \ref{case3a} and \ref{case3b}), which implies that $\dim \Ft=0$ so a three-step extension of $\Fn$ does not exist.\\

Now, assume that $\GG$ has inadmissible labels $z_{b_1},\dots,z_{b_n}$, and also assume that a three-step extension of $\Fn$ does not exist, i.e. $[v_i,z_{b_{\ell}}]=0 \; \forall v_i \in \Fv$ and $\forall \ell=1,\dots,n$. Now if we add an inadmissible label  $z_{b_{n+1}} \in C_{pos}$, we see that for any $v_1,v_2,v_3\in \Fv$,
$[v_1,[v_2,v_3]]+[v_2,[v_3,v_1]]+[v_3,[v_1,v_2]] =0$\\
$\implies ([v_1,\epsilon_{2,3}^{b_{n+1}} z_{b_{n+1}}]+[v_2,\epsilon_{3,1}^{b_{n+1}} z_{b_{n+1}}]+[v_3,\epsilon_{1,2}^{b_{n+1}} z_{b_{n+1}}])+\sum_{{\ell}=0}^n ([v_1,\epsilon_{2,3}^{b_{\ell}} z_{b_{\ell}}]+[v_2,\epsilon_{3,1}^{b_{\ell}} z_{b_{\ell}}]+[v_3,\epsilon_{1,2}^{b_{\ell}} z_{b_{\ell}}])=0$ \hfill by linearity of the bracket\\
$\implies [v_1,\epsilon_{2,3}^{b_{n+1}} z_{b_{n+1}}]+[v_2,\epsilon_{3,1}^{b_{n+1}} z_{b_{n+1}}]+[v_3,\epsilon_{1,2}^{b_{n+1}} z_{b_{n+1}}]=0$ \hfill by induction hypothesis\\
$\implies [v_i,z_{b_{n+1}}]=0 \; \forall v_i \in \Fv$ \hfill because the proof of the base case of this induction hypothesis showed that this is the result if the Jacobi identity is satisfied for any inadmissible label. \end{Proof}

\begin{Proof}[Proof of Thm. \ref{3stepthm} (nilpotency)] $[\widehat{\Fn},\widehat{\Fn}]=\Fz \oplus \Ft$ and $\widehat{\Fn}^{(2)}=[\widehat{\Fn},[\widehat{\Fn},\widehat{\Fn}]]=\Ft \subseteq Z(\widehat{\Fn})$. Therefore, $\widehat{\Fn}$ is a 3-step nilpotent Lie algebra. \end{Proof}


\section{Lie Algebras Associated with a Gassmann-Sunada Triple} \label{lags}
We continue the notation of \S\ref{defns} and \S\ref{23step}. The above construction does not require us to begin with a Gassmann-Sunada triple, but some interesting results occur when we look at the Lie algebras associated with a pair of Schreier graphs of a Gassmann-Sunada triple. Recall from Remark~\ref{onb} that in this paper, we will take the union of the set of vertices, the set of labels in $C_{pos}$, and the set $\{t_{r_{k,1}},\,t_{r_{k,2}}: \; z_{r_k} \mbox{ is admissible and } t_{r_{k,\ell}}\neq 0\}$ to be an orthonormal basis for $\widehat{\Fn}=\Fv\oplus\Fz\oplus\Ft$.

\begin{Proposition}\cite[Lecture 4]{CG} Let $(G,H_1,H_2)$ be a Gassmann-Sunada triple and $\GG_1$ and $\GG_2$ the pair of Schreier graphs associated with this triple. Let $\alpha_i$ be the group action of $G$ on $\Fv$ as in Remark~\ref{action}, which will be unitary under the assumed metric given in Remark \ref{onb}. Because $H_1$ and $H_2$ are almost conjugate subgroups of $G$, the representations $\alpha_1$ and $\alpha_2$ are unitarily equivalent, i.e. there exists a unitary operator $T: \Fv_1 \to \Fv_2$ such that $T(\alpha_1(x)(H_1 g))=\alpha_2(x)(T(H_1 g))$ for all $x \in G$ and for all $H_i g \in \Fv_i$, $i=1,2$. This operator $T$ is referred to as the transplantation or intertwining operator. For more information, see \cite{CG}. \end{Proposition}

Given a two-step nilpotent Lie algebra $\Fn=\Fv\oplus \Fz$, where $\Fv$ and $\Fz$ are inner product spaces, we can define an operator $j:\Fz \to so(\Fv)$ given by $j(z)(v)=(ad\; v)^* z$ where $(ad\; v)(w)=[v,w]$ and * denotes the adjoint operator with respect to the given inner product. In other words $j(z)v$ is the unique element in $\Fv$ such that
$$<j(z)v,w>=<z,[v,w]> \mbox{ for all } w \mbox{ in } \Fv$$ 
(See Eberlein \cite{eber}).\\

\begin{Theorem}\label{joper} The j-operator on the 2-step nilpotent metric Lie algebra, $\Fn=\Fv\oplus\Fz$, associated with a Schreier graph by Construction \ref{constr2} is given by, $\forall z \in \Fz \mbox{ and } \forall v \in \Fv$, $$j(z)v=\alpha(z)(v)-\alpha(z^{-1})(v).$$ \end{Theorem}
\begin{Proof} Fix basis elements $v \in \Fv$ and $z \in \Fz$.  Let $w$ be a basis element in $\Fv$. Then,
\begin{align*} <j(z)v,w> &= <z,[v,w]>\\
&= \left\{ \begin{array}{rcl} <z,z>=1 &,& \mbox{ if } w=\alpha(z)(v) \mbox{ and } \alpha(z)(v) \neq \alpha(z^{-1})(v)\\
<z,-z>=-1 &,& \mbox{ if } w=\alpha(z^{-1})(v) \mbox{ and } \alpha(z)(v) \neq \alpha(z^{-1})(v)\\
<z,z-z>=0 &,& \mbox{ if } w= \alpha(z)(v) = \alpha(z^{-1})(v)\\
0 &,& \mbox{ otherwise.} \end{array} \right. \end{align*}
Recall from Remark~\ref{possibilities} that this covers all cases that can occur on a Schreier graph.
On the other hand,\\
$<\alpha(z)(v)-\alpha(z^{-1})(v),w> = <\alpha(z)(v),w>-<\alpha(z^{-1})(v),w>$\\
$= \left\{ \begin{array}{rcl} 1-0=1 &,& \mbox{ if } w=\alpha(z)(v) \mbox{ and } \alpha(z)(v) \neq \alpha(z^{-1})(v)\\
0-1=-1 &,& \mbox{ if } w=\alpha(z^{-1})(v) \mbox{ and } \alpha(z)(v) \neq \alpha(z^{_1})(v)\\
<0,w>=0 &,& \mbox{ if } w= \alpha(z)(v) = \alpha(z^{-1})(v)\\
0 &,& \mbox{otherwise.} \end{array} \right.$\\
Since this is true for any basis elements $w$ in $\Fv$ and by the uniqueness and linearity of the inner product, this implies that $j(z)v=\alpha(z)(v)-\alpha(z^{-1})(v)$. \end{Proof}

\begin{Theorem} Starting with a pair of Schreier graphs coming from a Gassmann-Sunada triple, let $(\Fn_1,j_1)$ and $(\Fn_2,j_2)$ be the associated pair of two-step nilpotent Lie algebras determined by Construction \ref{constr2}, and let $T$ be the unitary intertwining operator guaranteed by the Gassmann-Sunada condition. Then, $$T(j_1(z)v)=j_2(z)(Tv) \; \forall z \in \Fz_1 \mbox{ and } \forall v \in \Fv_1.$$ \end{Theorem}
\begin{Proof} \begin{align*} T(j_1(z)v) &= T(\alpha_1(z)(v)-\alpha_1(z^{-1})(v))\\
&= T(\alpha_1(z)(v))-T(\alpha_1(z^{-1})(v))\\
&= \alpha_2(z)(Tv)-\alpha_2(z^{-1})(Tv)\\
&= j_2(z)(Tv) \tag*{\qedhere} \end{align*} \end{Proof}
\begin{Corollary} Starting with a pair of Schreier graphs coming from a Gassmann-Sunada triple, let $(\Fn_1,j_1)$ and $(\Fn_2,j_2)$ be the associated pair of two-step nilpotent metric Lie algebras determined by Construction \ref{constr2} with the metric defined in Remark \ref{onb}. Then, $(\Fn_1,j_1)$ is isometric to $(\Fn_2,j_2)$. \end{Corollary}
\begin{Proof} Using \cite[Lect. 8, Prop. 4.6]{CG}, we get $(\Fn_1,j_1)$ is isomorphic to $(\Fn_2,j_2)$ by $\widetilde{T}:=T\oplus Id$. 
Also $\forall v,v' \in \Fv_1 \mbox{ and } \forall z,z' \in \Fz_1$,\\ $<(v,z),(v',z')>_1=<v,v'>_1+<z,z'>_1$\\
$=<T(v),T(v')>_2+<Id(z),Id(z')>_2=<\widetilde{T}(v,z),\widetilde{T}(v',z')>_2$.  \end{Proof}

While the pair of two-step nilpotent Lie algebras associated with a Gassmann-Sunada triple are always isometric, the three-step nilpotent Lie algebra extensions determined by Construction \ref{constr3} need not be.
\begin{Theorem} The pair of three-step nilpotent Lie algebras given in Example~\ref{brooksex3} from \S\ref{23step} are non-isometric. \end{Theorem}
\begin{Proof} For the full proof, see the appendix below.  The idea of the proof is that we assume that there exists $\phi$ that is an isometry between $\widehat{\Fn_1}$ and $\widehat{\Fn_2}$. Then using the properties of Lie algebra isometries listed below, we obtain a contradiction.  Therefore, the two Lie algebras are non-isometric.
\begin{enumerate} \item $\phi:\Fv \to \Fv, \Fz \to \Fz, \mbox{ and } \Ft \to \Ft$.
\item $\phi$ has to preserve the ascending central series.
\item The columns (and rows) of the matrix $\phi$ must be orthonormal to each other.
\item $\phi([x,y]_1)=[\phi(x),\phi(y)]_2$ for all $x,y \in \widehat{\Fn_1}$. \qedhere
\end{enumerate} \end{Proof}
Note: Because there is a choice in constructing the 3-step nilpotent Lie algebra, a similar argument shows that the following variations on $\widehat{\Fn_2}$ are also non-isometric to $\widehat{\Fn_1}$:\\
\begin{enumerate} \item Interchanging $t$ and $-t$, i.e. $[v_3,z_r]=-t$ and $[v_6,z_r]=t$.
\item Switching the $t$ and $0$ components, i.e. $[v_3,z_r]=0,\, [v_5,z_r]=t,\, [v_6,z_r]=0, \mbox{ and } [v_7,z_r]=-t$.
\item Switching the $t$ and $0$ components and then interchanging $t$ and $-t$.
\end{enumerate}


\section{Appendix}\label{app}
Let $\widehat{\Fn_1}$ and $\widehat{\Fn_2}$ be the three step nilpotent Lie algebras given in Example~\ref{brooksex3}. Assume that $\phi:\;\widehat{\Fn_1}\to\widehat{\Fn_2}$ is an isometry, where the entries of the matrix $\phi$ with respect to the orthonormal basis $\{v_1,\dots,v_7,z_r,z_b,t\}$ for both $\widehat{\Fn_1}$ and $\widehat{\Fn_2}$ are $\left(\alpha_{i,j}\right)_{i,j=1}^{10}$.\\
We begin by computing the ascending central series of the two Lie algebras, obtaining the following:
\begin{align} Z(\widehat{\Fn_1}):=\{w\in \widehat{\Fn_1}:\, [w,\widehat{\Fn_1}]=0\} &=span_{\RR} \{v_1+v_2+\dots+v_6,v_7,z_b,t\}\label{zn1}\\ 
Z(\widehat{\Fn_2}):=\{w\in \widehat{\Fn_2}:\, [w,\widehat{\Fn_2}]=0\} &= span_{\RR} \{v_1+v_2+v_5+v_7,v_3+v_4+v_6,z_b,t\}\label{zn2}\\
Z_2(\widehat{\Fn_1}):=\{w\in \widehat{\Fn_1}:\, [w,\widehat{\Fn_1}]\subseteq Z(\widehat{\Fn_1})\} &=span_{\RR} \{v_1,v_2+v_4,v_3,v_5+v_6,v_7,z_r,z_b,t\}\label{z2n1}\\ 
Z_2(\widehat{\Fn_2}):=\{w\in \widehat{\Fn_2}:\, [w,\widehat{\Fn_2}]\subseteq Z(\widehat{\Fn_2})\} &=span_{\RR} \{v_1,v_2,v_3+v_6,v_4,v_5+v_7,z_r,z_b,t\}\label{z2n2}\end{align} 
Next, we use the assumption that $\phi$ is an isometry to obtain the following properties about the matrix $\phi$:
\begin{align} \phi \mbox{ an isometry} &\implies \mbox{ the matrix } \phi \mbox{ is orthonormal}\label{eqn2}\\
\phi:\Ft\to\Ft &\implies \alpha_{10,j}=\alpha_{j,10}=0 \mbox{ for } j=1,\dots,9\label{eqn1a}\\
\phi:\Fz\to\Fz &\implies \alpha_{8,j}=\alpha_{j,8}=\alpha_{9,j}=\alpha_{j,9}=0 \mbox{ for } j=1,\dots,7\label{eqn1b} \end{align} so that $\phi$ now is of the form
$$\left(
  \begin{array}{ccc|c|c}
    & & & & \\
		 &\mbox{\Huge A} & &0 & 0\\
     & & & &\\ \hline 
			& 0 & & \mbox{\large B} & 0\\ \hline
      & 0 & & 0 & \mbox{C}
  \end{array}
\right)$$
where A is of size $7\times7$, B is $2\times2$, and C is $1\times1$.\\
Finally, we use a combination of the above results along with the property that $\phi([x,y]_1)=[\phi(x),\phi(y)]_2$ for all $x,y \in \widehat{\Fn_1}$ to obtain relations about the various entries in $\phi$:
\begin{align} &\alpha_{10,10}=\pm1 \mbox{ (by \ref{eqn2} and \ref{eqn1a})}\label{eqnX}\\
\phi(z_b) \in Z(\widehat{\Fn_2}) &\implies \alpha_{8,9}=0 \mbox{ (by \ref{zn1} and \ref{zn2})}\label{eqn3a}\\
					&\implies \alpha_{9,9},\alpha_{8,8} \in \{\pm 1\} \mbox{ (by \ref{eqn2})}\label{eqn3b}\\
					&\mbox{ and } \alpha_{9,8}=0 \mbox{ (by \ref{eqn2})}\label{eqn3c}\\
\phi(v_7) \in Z(\widehat{\Fn_2}) &\implies \alpha_{1,7}=\alpha_{2,7}=\alpha_{5,7}=\alpha_{7,7}\notag\\ 
						&\mbox{ and }\hspace{.2cm} \alpha_{3,7}=\alpha_{4,7}=\alpha_{6,7} \mbox{ (by \ref{zn1} and \ref{zn2})} \label{eqn4}\\
\phi(v_1) \in Z_2(\widehat{\Fn_2}) &\implies \alpha_{3,1}=\alpha_{6,1} \mbox{ and } \alpha_{5,1}=\alpha_{7,1} \mbox{ (by \ref{z2n1} and \ref{z2n2})}\label{eqn6v1}\\
\phi(v_3) \in Z_2(\widehat{\Fn_2}) &\implies \alpha_{3,3}=\alpha_{6,3} \mbox{ and } \alpha_{5,3}=\alpha_{7,3} \mbox{ (by \ref{z2n1} and \ref{z2n2})}\label{eqn6v3}\\	
\phi(v_2+v_4) \in Z_2(\widehat{\Fn_2}) &\implies \alpha_{7,4}=\alpha_{5,2}+\alpha_{5,4}-\alpha_{7,2} \mbox{ (by \ref{z2n1} and \ref{z2n2})}\label{eqn6v2}\\
\phi(v_5+v_6) \in Z_2(\widehat{\Fn_2}) &\implies \alpha_{7,6}=\alpha_{5,5}+\alpha_{5,6}-\alpha_{7,5} \mbox{ (by \ref{z2n1} and \ref{z2n2})}\label{eqn6v5}\end{align} 
\begin{align} [\phi(v_2),\phi(z_r)]=\phi(t)=\alpha_{10,10}t &\implies \alpha_{6,2}=\alpha_{3,2}-\alpha_{8,8}\alpha_{10,10}\label{eqn7a}\\
[\phi(v_4),\phi(z_r)]=\phi(-t)=-\alpha_{10,10}t &\implies \alpha_{6,4}=\alpha_{3,4}+\alpha_{8,8}\alpha_{10,10}\label{eqn7b}\\
[\phi(v_5),\phi(z_r)]=\phi(0)=0 &\implies \alpha_{3,5}=\alpha_{6,5}\label{eqn8a}\\
[\phi(v_6),\phi(z_r)]=\phi(0)=0 &\implies \alpha_{3,6}=\alpha_{6,6}\label{eqn8b} \end{align}

\begin{align} (\mbox{row }3)\cdot(\mbox{row }6) &= 0 \mbox{ (by \ref{eqn2})} \notag\\*
&\implies \alpha_{3,4}=\alpha_{3,2}-\alpha_{8,8} \alpha_{10,10} \label{eqn9a}\\
&\implies \alpha_{6,4}=\alpha_{3,2} \mbox{ (by \ref{eqn4}, \ref{eqn6v1}, \ref{eqn6v3}, \ref{eqn7a}, \ref{eqn7b}, \ref{eqn8a}, \ref{eqn8b})} \label{eqn9b}\\
(\mbox{row }k)\cdot(\mbox{row }3) &= (\mbox{row }k)\cdot(\mbox{row }6) \mbox{ for $k=1,2,4,5,7$ (by \ref{eqn2})} \notag\\
&\implies \alpha_{k,4}=\alpha_{k,2} \mbox{ for } k=1,2,4,5,7 \mbox{ (by \ref{eqn4}, \ref{eqn6v1}, \ref{eqn6v3}, \ref{eqn7a}, \ref{eqn7b}, \ref{eqn8a}, \ref{eqn8b})} \label{eqn10a}\\
&\implies \alpha_{5,2}=\alpha_{5,4}=\alpha_{7,2}=\alpha_{7,4} \notag\\
&\implies \alpha_{7,4}=\alpha_{7,2}=2\alpha_{5,2}-\alpha_{7,2} \mbox{ (by \ref{eqn6v2})}  \label{eqn10b}\end{align}

\begin{align} [\phi(v_2),\phi(v_6)] &=\phi(-z_r)=-\alpha_{8,8}z_r, \mbox{ just looking at the $z_r$-coefficient} \notag\\
&\implies \sum_{i<j} (\alpha_{i,2}\alpha_{j,6}-\alpha_{j,2}\alpha_{i,6})[v_i,v_j]=-\alpha_{8,8}z_r \notag\\
&\implies -(\alpha_{3,2}\alpha_{5,6}-\alpha_{5,2}\alpha_{3,6})+(\alpha_{3,2}\alpha_{7,6}-\alpha_{7,2}\alpha_{3,6})\notag\\
&\hspace{1cm} -(\alpha_{5,2}\alpha_{6,6}-\alpha_{6,2}\alpha_{5,6})-(\alpha_{6,2}\alpha_{7,6}-\alpha_{7,2}\alpha_{6,6})=-\alpha_{8,8} \notag\\
&\implies \alpha_{7,6}=\alpha_{5,6}-\alpha_{10,10} \mbox{ (by \ref{eqn7a}, \ref{eqn8b}, \ref{eqn10b})} \label{eqn11a}\\
&\implies \alpha_{7,5}=\alpha_{5,5}+\alpha_{10,10} \mbox{ (by \ref{eqn6v5})} \label{eqn11b}\end{align}

\begin{align} [\phi(v_2),\phi(v_4)] &=\phi(0)=0, \mbox{ just looking at the $z_b$-coefficient} \notag\\
&\implies \sum_{i<j} (\alpha_{i,2}\alpha_{j,4}-\alpha_{j,2}\alpha_{i,4})[v_i,v_j]=0 \notag\\
&\implies (\alpha_{1,2}\alpha_{2,4}-\alpha_{2,2}\alpha_{1,4})-(\alpha_{1,2}\alpha_{5,4}-\alpha_{5,2}\alpha_{1,4})\notag\\
&\hspace{1cm} +(\alpha_{2,2}\alpha_{5,4}-\alpha_{5,2}\alpha_{2,4})-(\alpha_{3,2}\alpha_{4,4}-\alpha_{4,2}\alpha_{3,4}) \notag\\
&\hspace{1cm} +(\alpha_{3,2}\alpha_{6,4}-\alpha_{6,2}\alpha_{3,4})-(\alpha_{4,2}\alpha_{6,4}-\alpha_{6,2}\alpha_{4,4})=0 \notag\\
&\implies \alpha_{4,2}=\alpha_{3,2}-1/2\alpha_{8,8} \alpha_{10,10} \mbox{ (by \ref{eqn7a}, \ref{eqn9a}, \ref{eqn9b}, \ref{eqn10a})} \label{eqn12} \end{align}

\begin{align} (\mbox{row }k)\cdot(\mbox{row }5) &= (\mbox{row }k)\cdot(\mbox{row }7) \mbox{ for $k=1,2,3,4,6$ (by \ref{eqn2})} \notag\\
&\implies \alpha_{k,5}=\alpha_{k,6} \mbox{ for $k=1,2,3,4,6$ (by \ref{eqn4}, \ref{eqn6v1}, \ref{eqn6v3}, \ref{eqn10b}, \ref{eqn11a}, \ref{eqn11b})} \label{eqn13}\\
||\mbox{row }5||-1 &= (\mbox{row }5)\cdot(\mbox{row }7) \mbox{ (by \ref{eqn2})} \notag\\
&\implies \alpha_{5,6}=\alpha_{5,5}+\alpha_{10,10} \mbox{ (by \ref{eqn4}, \ref{eqn6v1}, \ref{eqn6v3}, \ref{eqn10b}, \ref{eqn11a}, \ref{eqn11b})} \label{eqn14a}\\
&\implies \alpha_{7,6} = \alpha_{5,5} \mbox{ (by \ref{eqn11a})} \label{eqn14b}\end{align}

\begin{align} [\phi(v_5),\phi(v_6)] &=\phi(0)=0, \mbox{ just looking at the $z_b$-coefficient} \notag\\
&\implies \sum_{i<j} (\alpha_{i,5}\alpha_{j,6}-\alpha_{j,5}\alpha_{i,6})[v_i,v_j]=0 \notag\\
&\implies (\alpha_{1,5}\alpha_{2,6}-\alpha_{2,5}\alpha_{1,6})-(\alpha_{1,5}\alpha_{5,6}-\alpha_{5,5}\alpha_{1,6})\notag\\
&\hspace{1cm} +(\alpha_{2,5}\alpha_{5,6}-\alpha_{5,5}\alpha_{2,6})-(\alpha_{3,5}\alpha_{4,6}-\alpha_{4,5}\alpha_{3,6}) \notag\\
&\hspace{1cm} +(\alpha_{3,5}\alpha_{6,6}-\alpha_{6,5}\alpha_{3,6})-(\alpha_{4,5}\alpha_{6,6}-\alpha_{6,5}\alpha_{4,6})=0 \notag\\
&\implies \alpha_{1,5}=\alpha_{2,5} \mbox{ (by \ref{eqn8a}, \ref{eqn8b}, \ref{eqn13}, \ref{eqn14a})} \label{eqn15} \end{align}

\begin{align} [\phi(v_4),\phi(v_5)] &=\phi(-z_r)=-\alpha_{8,8}z_r, \mbox{ just looking at the $z_b$-coefficient} \notag\\
&\implies \sum_{i<j} (\alpha_{i,4}\alpha_{j,5}-\alpha_{j,4}\alpha_{i,5})[v_i,v_j]=-\alpha_{8,8}z_r \notag\\
&\implies (\alpha_{1,4}\alpha_{2,5}-\alpha_{2,4}\alpha_{1,5})-(\alpha_{1,4}\alpha_{5,5}-\alpha_{5,4}\alpha_{1,5})\notag\\
&\hspace{1cm} +(\alpha_{2,4}\alpha_{5,5}-\alpha_{5,4}\alpha_{2,5})-(\alpha_{3,4}\alpha_{4,5}-\alpha_{4,4}\alpha_{3,5}) \notag\\
&\hspace{1cm} +(\alpha_{3,4}\alpha_{6,5}-\alpha_{6,4}\alpha_{3,5})-(\alpha_{4,4}\alpha_{6,5}-\alpha_{6,4}\alpha_{4,5})=0 \notag\\
&\implies \alpha_{1,2}\alpha_{1,5}-\alpha_{2,2}\alpha_{1,5}-\alpha_{1,2}\alpha_{5,5}+\alpha_{2,2}\alpha_{5,5} \notag\\
&\hspace{1cm} =-\alpha_{4,5}\alpha_{8,8}\alpha_{10,10}+\alpha_{3,5}\alpha_{8,8}\alpha_{10,10} \mbox{ (by \ref{eqn7b}, \ref{eqn8a}, \ref{eqn10a}, \ref{eqn15})} \label{eqn17a} \end{align}

\begin{align} [\phi(v_4),\phi(v_6)] &=\phi(z_r+z_b)=\alpha_{8,8}z_r+\alpha_{9,9}z_b, \mbox{ just looking at the $z_b$-coefficient} \notag\\
&\implies \sum_{i<j} (\alpha_{i,4}\alpha_{j,6}-\alpha_{j,4}\alpha_{i,6})[v_i,v_j]=-\alpha_{8,8}z_r+\alpha_{9,9}z_b \notag\\
&\implies (\alpha_{1,4}\alpha_{2,6}-\alpha_{2,4}\alpha_{1,6})-(\alpha_{1,4}\alpha_{5,6}-\alpha_{5,4}\alpha_{1,6})\notag\\
&\hspace{1cm} +(\alpha_{2,4}\alpha_{5,6}-\alpha_{5,4}\alpha_{2,6})-(\alpha_{3,4}\alpha_{4,6}-\alpha_{4,4}\alpha_{3,6}) \notag\\
&\hspace{1cm} +(\alpha_{3,4}\alpha_{6,6}-\alpha_{6,4}\alpha_{3,6})-(\alpha_{4,4}\alpha_{6,6}-\alpha_{6,4}\alpha_{4,6})=\alpha_{9,9} \notag \\
&\implies \alpha_{1,2}\alpha_{1,5}-\alpha_{2,2}\alpha_{1,5}-\alpha_{1,2}\alpha_{5,5}+\alpha_{2,2}\alpha_{5,5}\notag\\
&\hspace{1cm} =\alpha_{1,2}\alpha_{10,10}-\alpha_{2,2}\alpha_{10,10}-\alpha_{4,5}\alpha_{8,8}\alpha_{10,10}+\alpha_{3,5}\alpha_{8,8}\alpha_{10,10}+\alpha_{9,9} \label{eqn17b}\\
&\hspace{1cm} \mbox{(by \ref{eqn7a}, \ref{eqn8a}, \ref{eqn9a}, \ref{eqn10a}, \ref{eqn13}, \ref{eqn14a})} \notag \end{align}

\begin{align} [\phi(v_2),\phi(v_5)] &=\phi(z_r-z_b)=\alpha_{8,8}z_r-\alpha_{9,9}z_b, \mbox{ just looking at the $z_b$-coefficient} \notag\\
&\implies \sum_{i<j} (\alpha_{i,2}\alpha_{j,5}-\alpha_{j,2}\alpha_{i,5})[v_i,v_j]=-\alpha_{8,8}z_r-\alpha_{9,9}z_b \notag\\
&\implies (\alpha_{1,2}\alpha_{2,5}-\alpha_{2,2}\alpha_{1,5})-(\alpha_{1,2}\alpha_{5,5}-\alpha_{5,2}\alpha_{1,5})\notag\\
&\hspace{1cm} +(\alpha_{2,2}\alpha_{5,5}-\alpha_{5,2}\alpha_{2,5})-(\alpha_{3,2}\alpha_{4,5}-\alpha_{4,2}\alpha_{3,5}) \notag\\
&\hspace{1cm} +(\alpha_{3,2}\alpha_{6,5}-\alpha_{6,2}\alpha_{3,5})-(\alpha_{4,2}\alpha_{6,5}-\alpha_{6,2}\alpha_{4,5})=-\alpha_{9,9} \notag \\
&\implies \alpha_{1,2}\alpha_{1,5}-\alpha_{2,2}\alpha_{1,5}-\alpha_{1,2}\alpha_{5,5}+\alpha_{2,2}\alpha_{5,5} \notag\\
&\hspace{1cm}=-\alpha_{3,5}\alpha_{8,8}\alpha_{10,10}+\alpha_{4,5}\alpha_{8,8}\alpha_{10,10}-\alpha_{9,9} \mbox{ (by \ref{eqn7a}, \ref{eqn8a}, \ref{eqn15})} \label{eqn17c} \end{align}

\begin{align} \mbox{So, } \alpha_{2,2} &=\alpha_{1,2}+\alpha_{9,9}\alpha_{10,10} \mbox{ (by \ref{eqn17a}, \ref{eqn17b})} \label{eqn17ab}\\
\mbox{and } \alpha_{4,5} &=\alpha_{3,5}+1/2\alpha_{8,8}\alpha_{9,9}\alpha_{10,10}\mbox{ (by \ref{eqn17a}, \ref{eqn17c})}\label{eqn17ac} \end{align}

\begin{align} [\phi(v_1),\phi(v_5)] &=\phi(z_b)=\alpha_{9,9}z_b, \mbox{ just looking at the $z_b$-coefficient} \notag\\*
&\implies \sum_{i<j} (\alpha_{i,1}\alpha_{j,5}-\alpha_{j,1}\alpha_{i,5})[v_i,v_j]=\alpha_{9,9}z_b \notag\\
&\implies (\alpha_{1,1}\alpha_{2,5}-\alpha_{2,1}\alpha_{1,5})-(\alpha_{1,1}\alpha_{5,5}-\alpha_{5,1}\alpha_{1,5})\notag\\
&\hspace{1cm} +(\alpha_{2,1}\alpha_{5,5}-\alpha_{5,1}\alpha_{2,5})-(\alpha_{3,1}\alpha_{4,5}-\alpha_{4,1}\alpha_{3,5}) \notag\\
&\hspace{1cm} +(\alpha_{3,1}\alpha_{6,5}-\alpha_{6,1}\alpha_{3,5})-(\alpha_{4,1}\alpha_{6,5}-\alpha_{6,1}\alpha_{4,5})=\alpha_{9,9} \notag \\
&\implies \alpha_{1,1}\alpha_{1,5}-\alpha_{2,1}\alpha_{1,5}-\alpha_{1,1}\alpha_{5,5}+\alpha_{2,1}\alpha_{5,5}\notag\\
&\hspace{1cm} =\alpha_{9,9} \mbox{ (by \ref{eqn6v1}, \ref{eqn8a}, \ref{eqn15})} \label{eqn1,5b} \end{align}

\begin{align} [\phi(v_1),\phi(v_6)] &=\phi(0)=0, \mbox{ just looking at the $z_b$-coefficient} \notag\\
&\implies \sum_{i<j} (\alpha_{i,1}\alpha_{j,6}-\alpha_{j,1}\alpha_{i,6})[v_i,v_j]=0 \notag\\
&\implies (\alpha_{1,1}\alpha_{2,6}-\alpha_{2,1}\alpha_{1,6})-(\alpha_{1,1}\alpha_{5,6}-\alpha_{5,1}\alpha_{1,6})\notag\\
&\hspace{1cm} +(\alpha_{2,1}\alpha_{5,6}-\alpha_{5,1}\alpha_{2,6})-(\alpha_{3,1}\alpha_{4,6}-\alpha_{4,1}\alpha_{3,6}) \notag\\
&\hspace{1cm} +(\alpha_{3,1}\alpha_{6,6}-\alpha_{6,1}\alpha_{3,6})-(\alpha_{4,1}\alpha_{6,6}-\alpha_{6,1}\alpha_{4,6})=0 \notag \\
&\implies \alpha_{1,1}\alpha_{1,5}-\alpha_{2,1}\alpha_{1,5}-\alpha_{1,1}\alpha_{5,5}+\alpha_{2,1}\alpha_{5,5}\notag\\
&\hspace{1cm}=\alpha_{1,1}\alpha_{10,10}-\alpha_{2,1}\alpha_{10,10} \mbox{ (by \ref{eqn6v1}, \ref{eqn8b}, \ref{eqn13}, \ref{eqn14a}, \ref{eqn15})} \label{eqn1,6b}\\
\mbox{So, } \alpha_{2,1} &= \alpha_{1,1}-\alpha_{9,9}\alpha_{10,10} \mbox{ (by \ref{eqn1,5b}, \ref{eqn1,6b})} \label{eqn18a}\\
\mbox{and } \alpha_{5,5} &= \alpha_{1,5}-\alpha_{10,10} \mbox{ (by \ref{eqn1,5b}, \ref{eqn1,6b}, \ref{eqn18a})} \label{eqn18b}\end{align}

\begin{align} [\phi(v_2),\phi(v_5)] &=\phi(z_r-z_b)=\alpha_{8,8}z_r-\alpha_{9,9}z_b, \mbox{ just looking at the $z_b$-coefficient} \notag\\
&\implies \sum_{i<j} (\alpha_{i,2}\alpha_{j,5}-\alpha_{j,2}\alpha_{i,5})[v_i,v_j]=-\alpha_{8,8}z_r-\alpha_{9,9}z_b \notag\\
&\implies (\alpha_{1,2}\alpha_{2,5}-\alpha_{2,2}\alpha_{1,5})-(\alpha_{1,2}\alpha_{5,5}-\alpha_{5,2}\alpha_{1,5})\notag\\
&\hspace{1cm} +(\alpha_{2,2}\alpha_{5,5}-\alpha_{5,2}\alpha_{2,5})-(\alpha_{3,2}\alpha_{4,5}-\alpha_{4,2}\alpha_{3,5}) \notag\\
&\hspace{1cm} +(\alpha_{3,2}\alpha_{6,5}-\alpha_{6,2}\alpha_{3,5})-(\alpha_{4,2}\alpha_{6,5}-\alpha_{6,2}\alpha_{4,5})=-\alpha_{9,9} \notag \\
&\implies -3/2\alpha_{9,9}=-\alpha_{9,9} \mbox{ (by \ref{eqn7a}, \ref{eqn8a}, \ref{eqn15}, \ref{eqn17ab}, \ref{eqn17ac}, \ref{eqn18b})} \notag\\
&\implies \alpha_{9,9}=0 \mbox{ which contradicts equation \ref{eqn3b}} \notag\end{align}
Therefore, $\widehat{\Fn_1}$ is not isometric to $\widehat{\Fn_2}$.       


\end{document}